\documentclass[a4paper, 10pt, reqno]{amsart}

\usepackage{ amssymb, amsmath, amsthm, textcomp}
\usepackage{graphicx, psfrag, setspace, subfigure, gensymb}
\usepackage{amsfonts}
\usepackage{bbm}
\usepackage{fix-cm}
\usepackage{comment}
\usepackage{tikz,tikz-cd}
\usetikzlibrary{matrix,arrows,decorations.pathmorphing,decorations.pathreplacing, calc}
\tikzset{commutative diagrams/diagrams={baseline=-2.5pt},commutative diagrams/arrow style=tikz}
\usepackage[colorlinks]{hyperref}
\usepackage{float, enumitem}
\usepackage{setspace} % needed for redefined \marginpar

\newcommand\Z{\mathbb Z}
\newcommand\C{\mathbb C}
\newcommand\R{\mathbb R}
\newcommand\N{\mathbb N}

\newcommand\A{\mathbb{A}}

\newcommand{\cC}{\mathcal{C}}

\newcommand{\cO}{\mathcal{O}}

\newcommand{\cS}{\mathcal{S}}
\newcommand{\cT}{\mathcal T}

\newcommand{\mfF}{\mathfrak{F}}

\newcommand\isoto{\stackrel{\sim}{\To}}
\newcommand\id{\mathrm 1}

\newcommand\into{\hookrightarrow}
\newcommand\To{\longrightarrow}

\renewcommand\Im{\operatorname{Im}}
\newcommand\Ker{\operatorname{Ker}}

\renewcommand\P{\mathbb P}

\newcommand\rk{\operatorname{rank}}
\newcommand\rank{\operatorname{rank}}

\newcommand\actson{\curvearrowright}

\newcommand{\al}[1]{\begin{align*}#1\end{align*}}

\newcommand{\beq}[1]{\begin{equation}\label{#1} }
\newcommand{\eeq}{\end{equation}}
\newcommand{\pgap}{\vspace{5pt}}

\theoremstyle{plain}
\newtheorem{prop}[equation]{Proposition}
\newtheorem{thm}[equation]{Theorem}
\newtheorem{lem}[equation]{Lemma}
\newtheorem{cor}[equation]{Corollary}
\newtheorem{conj}[equation]{Conjecture}

\theoremstyle{remark}
\newtheorem{rem}[equation]{Remark}

\theoremstyle{definition}
\newtheorem{defn}[equation]{Definition}

\newtheorem{defnlem}[equation]{Definition/Lemma}

\newtheorem{examplex}[equation]{Example}
\newenvironment{eg}
  {\pushQED{\qed}\examplex}
  {\popQED\endexamplex}

\makeatletter \@addtoreset{equation}{section} \makeatother

\newtheorem{keythm}{Theorem}

\newtheorem{keyconj}[keythm]{Conjecture}

\let\oldtocsection=\tocsection
\let\oldtocsubsection=\tocsubsection
\let\oldtocsubsubsection=\tocsubsubsection
\renewcommand{\tocsection}[3]{\hspace{0em}\oldtocsection{#1}{#2}{#3}}
\renewcommand{\tocsubsection}[3]{ \hspace{1em} \oldtocsubsection{#1}{\small{#2}}{\small{#3}} }
\renewcommand{\tocsubsubsection}[3]{\hspace{2em}\oldtocsubsubsection{#1}{\small{#2}}{\small{#3}}}

\newcommand{\marginparstretch}{0.6}
\let\oldmarginpar\marginpar
\renewcommand\marginpar[1]{\-\oldmarginpar[\framebox{\setstretch{\marginparstretch}\begin{minipage}{\marginparwidth}{\raggedleft\scriptsize #1}\end{minipage}}]{\framebox{\setstretch{\marginparstretch}\begin{minipage}{\marginparwidth}{\raggedright\scriptsize #1}\end{minipage}}}}

\newcommand\SOD[1]{\big\langle\,#1\,\big\rangle}

\begin{document}

\title{Discriminants and semi-orthogonal decompositions}
\author{Alex Kite and Ed Segal} 

\maketitle

%\emph{Communicating author:}

%\emph{Ed Segal, University College London, \href{mailto:e.segal@ucl.ac.uk}{e.segal@ucl.ac.uk}}

\begin{abstract} The derived categories of toric varieties admit semi-orthogonal decompositions coming from wall-crossing in GIT. We prove that these decompositions satisfy a Jordan-H\"older property: the subcategories that appear, and their multiplicities, are independent of the choices made.

For Calabi-Yau toric varieties wall-crossing instead gives derived equivalences and autoequivalences, and mirror symmetry relates these to monodromy around the GKZ discriminant locus. We formulate a conjecture equating intersection multiplicities in the discriminant with the multiplicities appearing in certain semi-orthogonal decompositions. We then prove this conjecture in some cases. 
\end{abstract}

\tableofcontents

\section{Introduction}

Let $X$ be a toric variety, constructed as a GIT quotient of a vector space $V$ by a torus $T$. There is a well-established theory \cite{Kaw, Seg, BFK, HL} that tells us how to produce semi-orthogonal decompositions of the derived category $D^b(X)$. We do it by considering other birational models of $X$, \emph{i.e.} crossing walls in the GIT problem $T\actson V$. If we cross to a quotient $X'$, and $K_{X'}$ is `more negative' than $K_X$, then $D^b(X)$ decomposes as
\beq{eq.SOD}D^b(X) = \SOD{D^b(X'), \;D^b(Z),\; ..., \;D^b(Z)}\eeq
where $Z$ is another toric variety of smaller dimension. We do this repeatedly until we arrive at a `minimal' chamber. Since the extra pieces are always equivalent to the derived category of a toric variety they themselves can be decomposed by the same procedure, and we get a recursive algorithm which terminates after a finite number of steps.

If $X$ is projective then the result of this algorithm is a full exceptional collection for $X$, \emph{i.e.} every piece of the final decomposition is equivalent to $D^b(\C)$. But for quasi-projective varieties there will usually be many different categories occuring, each one with some multiplicity. Moreover the decomposition is not unique; at each step of the algorithm one may have a choice about which wall to cross through and these choices result in different decompositions. The main technical result of this paper is the following Jordan-H\"older type theorem:

\begin{keythm}[Theorem \ref{thm.JH}] Let $X$ be a toric variety. If we decompose $D^b(X)$ using the wall-crossing algorithm then the subcategories occuring in the final decomposition, and their multiplicities, are independent of all choices.
\end{keythm}

This result is not particularly hard to prove and neither is it an abstract result; we prove it by analysing the algorithm.  But it is notable that the Jordan-H\"older property does not hold for semi-orthogonal decompositions in general \cite{BBS, Kuz}. 
\pgap

Our real motivation for proving the theorem above was to be able to understand a conjecture appearing in a physics paper by Aspinwall--Plesser--Wang \cite{APW}.  Part of what they state is already understood in the mathematical literature but there remains a significant unsolved problem which we are able to formulate precisely using our theorem (Conjecture \ref{conj.mainconj}). This generalizes a conjecture made by Halpern-Leistner--Shipman \cite{HLSh}. 

 We will use the remainder of this introduction to explain the motivation and context for this conjecture

\subsection{Spherical functors from wall-crossing}\label{sec.sphericalintro}

Our conjecture concerns the special case when the torus action $T\actson V$ is through the subgroup $SL(V)$. In this case all the GIT quotients $X$ will be Calabi-Yau, meaning $K_X\cong \cO_X$, and not projective. In this situation the wall-crossing theory does not provide any decompositions of $D^b(X)$, instead it proves that all the GIT quotients are derived equivalent since the decomposition \eqref{eq.SOD} just becomes $D^b(X)=D^b(X')$. However the category $D^b(Z)$ still has an important role.

The derived equivalence between $X$ and $X'$ is not unique, the theory gives us multiple equivalences for every wall-crossing, and by composing them we get autoequivalences of $D^b(X)$. From work of Halpern-Leistner--Shipman \cite{HLSh} it is known that each of these autoequivalences can be described as a twist $T_F$ around a spherical functor
$$F: D^b(Z) \to D^b(X)$$
where $Z$ is the same toric variety that appears in \eqref{eq.SOD}. 

 By combining these, and the Picard groups of each GIT quotient, we can get many autoequivalences of $D^b(X)$.  So the interesting problem becomes to understand this large group of autoequivalences.

\subsection{FI parameter spaces}\label{sec.FIPSintro}

Now we explain some heuristics from physics and mirror symmetry. In string theory the data of $T$ acting on $V$ determines an abelian gauged linear sigma model, a widely studied class of $N=(2,2)$ superconformal field theories. In this theory there are certain important parameters called complexified Fayet-Iliopoulos parameters, they take values in a complex manifold which we call the FI Parameter Space (FIPS). They are related to stability conditions in the GIT problem and in certain limiting regions of the FIPS the theory reduces to a sigma model whose target is one of the quotients $X$. In physical terminology $X$ is a \emph{phase} of the model. In this region we can identify the FI parameters with the complexified K\"ahler moduli of $X$ so the FIPS is closely related to the \emph{extended} or \emph{stringy} K\"ahler moduli space of $X$.\footnote{The FIPS is not quite the same as the SKMS, the latter should be intrinsic to $X$ whereas the former depends on its presentation as $V/\!/T$. Also note that the SKMS is expected to be a complex submanifold of the space of Bridgeland stability conditions; on the mirror side this is the difference between small and big quantum cohomology. The FIPS is easier to compute than either the SKMS or the space of stability conditions.} 

Under mirror symmetry the FI parameters become complex parameters, so the FIPS is the base of the mirror family. Since toric mirror symmetry has a mathematically precise formulation this gives us a rigorous definition of the FIPS: it's the complement of the GKZ discriminant locus $\nabla$ inside the dual torus $T^\vee$ (Section \ref{sec.discriminants}). It is helpful to think of  $T^\vee$ as an open subset of the secondary toric variety $\mfF$ and to take the closure $\overline{\nabla}\subset \mfF$, because then the phases correspond to the toric fixed points in $\mfF$. From this point-of-view the FIPS is obtained by deleting $\overline{\nabla}$ and the toric boundary from $\mfF$. 

The mirror family is a locally-trivial family of symplectic manifolds over the FIPS with fibre $\check{X}$. The monodromy of this family gives an action of $\pi_1(FIPS)$ on $\check{X}$ as symplectomorphisms, and hence as autoequivalences of the Fukaya category $\mathrm{Fuk}(\check{X})$. On the mirror side this predicts an action:
$$\pi_1(FIPS) \actson D^b(X)$$
This is the `B-brane monodromy'.  Examples and physical calculations suggest that this is essentially the group of autoequivalences that arise via wall-crossing as described in Section \ref{sec.sphericalintro}.  This prediction appears in many places in the maths and physics literature (\emph{e.g.} \cite{HHP, HW, HLSam, HLSh}) and has been verified for some examples \cite{DS, Kite}. It seems to be a difficult problem to verify it in general, mainly because it is hard to understand $\pi_1(FIPS)$.

\subsection{The rank 1 case}\label{sec.rank1} 

The case where $T=\C^*$ is quite well-known and easy to understand directly. In this case there are two possible phases which we denote by $X_\pm$. If we split $V$ by weights as $V_+\oplus V_0\oplus V_-$ then it's easy to see that $X_\pm$ is a vector bundle over $\P V_\pm \times V_0$, where $\P V_\pm$ is a weighted projective space.

In this rank 1 case the discriminant locus is always a single point $\delta$ so the FIPS is $\C^*\setminus \delta$ (see Example \ref{eg.rank1discriminant}). Or we can say that the secondary toric variety $\mfF$ is a $\P^1$ and that the FIPS is obtained from it by deleting the two toric fixed points and one more non-fixed point. The phase $X_+$ corresponds to the region near one of the toric fixed points, and the loop around that fixed point simply acts as $\otimes \cO(1)$ on $D^b(X_+)$.

  More interesting is the loop around the non-fixed point $\delta$ - often called the \emph{conifold point} -  which  corresponds to wall-crossing to $X_-$ and back again. If there are no zero weights then the resulting autoequivalence is the twist $T_S$ around a spherical object
$$S = \cO_{P V_+}$$
given by the sky-scraper sheaf along the zero section in $X_+$. If there are zero weights we upgrade this to a twist around the spherical functor
$$F: D^b(V_0) \to D^b(X_+)$$
given by pulling-up to $ \P V_+\times V_0$ and then pushing-forward along the inclusion into $X_+$. In the notation of Section \ref{sec.sphericalintro} the variety $Z$ is $V_0$. 

%In this case there is nothing to prove since $\pi(FIPS)$ is free. But it is more satisfying to consider a fundamental groupoid of FIPS with two basepoints, one for each phase, then the thing to check is that the wall-crossing equivalences interact in the correct way with the Picard groups of each phase.

\begin{rem}\label{rem.orbifoldpi1}If there is only one positive weight then $X_+$ is an affine orbifold and $\mathrm{Pic}(X_+)$ is a finite cyclic group $\Z/k$. In this case it's sensible to allow that toric fixed point as part of the FIPS. The reason is that $\mfF$ is (if we're careful) an orbifold $\P^1$ and this fixed point has isotropy group $\Z/k$, so we get an action of the orbifold fundamental group. 

 This subtlety is interesting in the rank 1 case since it is occurs in the well-known `Calabi-Yau/Landau-Ginsburg correspondence'. In higher rank it happens very rarely and is of no significance for this paper. For us the FIPS will contain none of the toric boundary and hence we can ignore any orbifold structure on $\mfF$. 
\end{rem}

\subsection{Components of the discriminant}\label{sec.componentsintro}

Suppose we have a higher rank torus $T\cong (\C^*)^r$. The discrimant locus $\nabla$ is now some hypersurface in $(\C^*)^r$ and it is usually the union of several irreducible components:
 $$\nabla = \nabla_0\cup... \cup \nabla_k$$
Aspinwall--Plesser--Wang \cite{APW} observed that there is a correspondence between these components $\nabla_i$ and certain toric varieties $Z_i$, built from subsets of the original toric data. They conjecture that for each phase $X$ there should be a spherical functor
\beq{eq.Fi}F_i: D^b(Z_i) \to D^b(X)\eeq
and that $T_{F_i}$ corresponds to the monodromy around the component $\nabla_i$. There is some deliberate ambiguity here; there is no canonical loop around $\nabla_i$ (even up to homotopy), so the functors $F_i$ are at best defined up to composition by autoequivalences.

\subsection{Factorizations and multiplicities}\label{sec.multsintro}

To understand this conjecture of \cite{APW} more clearly we pick two adjacent chambers of the secondary fan, separated by a wall $W$. This is the situation we discussed in Section \ref{sec.sphericalintro}. The two chambers give two phases $X_\pm$ which are derived equivalent, and we get an autoequivalence of $D^b(X_+)$ which is the twist around a spherical functor 
\beq{eq.FW}F: D^b(Z) \to D^b(X_+)\eeq
 for some smaller toric variety $Z$. 

In the secondary toric variety $\mfF$ our wall $W$ corresponds to a rational curve $C_W$ connecting the toric fixed points  corresponding to our two phases.  It turns out that the discriminant locus $\overline{\nabla}$ always intersects $C_W$ in a single point $\delta$ (Corollary \ref{cor.onepointinCW}). This is the same picture that we saw in Section \ref{sec.rank1}, and the reason for this is that by focusing on a single wall-crossing we are essentially reducing to a rank 1 GIT problem. There is a the 1-parameter subgroup $\lambda_W\subset T$ normal to the wall and it is only stability with respect to $\lambda_W$ that is changing. So, just as in the rank 1 case, a loop in $C_W$ that goes around the point $\delta$ should correspond to the autoequivalence $T_F$.

However, $C_W$ is not part of the FIPS since it lies in the toric boundary of $\mfF$. To get an actual element of $\pi_1(FIPS)$ we have to perturb $C_W$ (or an open subset of it) off the toric boundary, and take a loop in the perturbed curve. 

\begin{figure}[!tbp]
  \centering
  \begin{minipage}[b]{0.5\textwidth}
    \includegraphics[width=\textwidth]{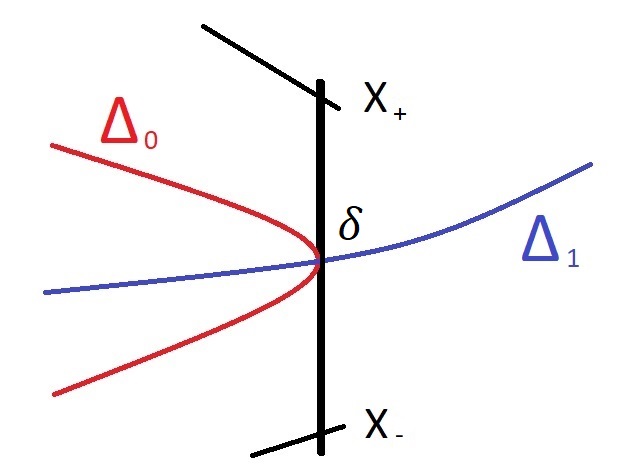}
   % \caption{(a)}
  \end{minipage}
  \hfill
  \begin{minipage}[b]{0.35\textwidth}
    \includegraphics[width=\textwidth]{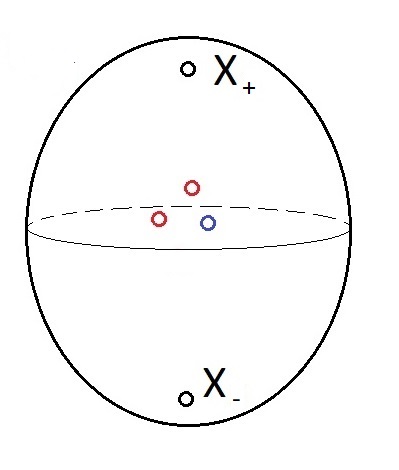}
%\caption{(b)}
   
  \end{minipage}
\caption{(L) A real picture of $C_W$ as the straight line connecting the two points marked by $X_\pm$. (R) A complex picture of a 2-sphere near to the rational curve $C_W$, where the point $\delta$ has split into three. A loop from $X_+$ to $X_-$ and back again will factor into two loops around $\Delta_0$ and one loop around $\Delta_1$.}
\label{fig.FIPS}
\end{figure}

When we do this peturbation the point $\delta$ may split into several points because $\overline{\nabla}$ typically meets $C_W$ with some multiplicity. This means that our element of $\pi_1(FIPS)$ is naturally a composite of several loops, one around each of our new missing points.  In fact each component $\overline{\nabla}_i$ might meet $C_W$ with multiplicity, and we can group the new missing points according to these components (see Figure \ref{fig.FIPS}). 

So the loop around $\delta$ naturally factors into several loops around the different components of $\nabla$, with each component possibly appearing multiple times.  This suggests that we should look for a corresponding factorization of the autoequivalence $T_F$. 
\pgap

This factorization does indeed exist. The toric variety $Z$ is not usually a Calabi-Yau, which means that $D^b(Z)$ (unlike $D^b(X)$) can be decomposed using the wall-crossing algorithm. Moreover, the subcategories that appear in this decomposition are always equivalent to $D^b(Z_i)$ where $Z_i$ is one of the varieties considered by Aspinwall--Plesser--Wang (Section \ref{sec.componentsintro}). So we get a semi-orthogonal decomposition
\beq{eq.SODW}D^b(Z) = \SOD{ D^b(Z_0), \; D^b(Z_0), \;  ..., \; D^b(Z_k), \; D^b(Z_k) } \eeq
where each $D^b(Z_i)$ occurs some number of times (possibly zero). The order of the factors here depends on the choices made in the algorithm, but by our Theorem \ref{thm.JH} the multiplicities do not. 

Halpern-Leistner--Shipman \cite{HLSh} observed that this decomposition gives us a factorization of the autoequivalence $T_F$. If we restrict the spherical functor $F$ \eqref{eq.FW} to each piece of $D^b(Z)$ then we again get a spherical functor, and $T_F$ is the composition of all the corresponding twists. This provides the spherical functors $F_i$ required by Aspinwall--Plesser--Wang  and matches with our discussion of loops in the FIPS.

However, for this story to make sense there is one essential numerical condition:

\begin{keyconj}[Conjecture \ref{conj.mainconj}]\label{conj.key} The multiplicity of $D^b(Z_i)$ in the decomposition \eqref{eq.SODW} agrees with the intersection multiplicity of $\overline{\nabla}_i$ with $C_W$. 
\end{keyconj}

We finish by proving our conjecture in some special cases, the strongest of which is:

\begin{keythm}[Theorem \ref{thm.rank2}] If the torus $T$ has rank 2 then Conjecture \ref{conj.key} holds.
\end{keythm}

\begin{rem}
A significant part of this story was already understood by Halpern-Leistner--Shipman. They only consider the case when $Z$ is projective, meaning that the decomposition of $D^b(Z)$ is actually a full exceptional collection, and they conjecture that the number of exceptional objects agrees with the intersection multiplicity of $\overline{\nabla}$ with $C_W$ \cite[Remark 4.7]{HLSh}. Our conjecture is a synthesis of theirs with the work of \cite{APW}. 
\end{rem}

\subsection{Acknowledgements}

E.S. would like to thank Paul Aspinwall and Lars Louder for helpful conversations. 
\pgap

This project has received funding from the European Research Council (ERC) under the European Union Horizon 2020 research and innovation programme (grant agreement No.725010). A.K. was supported by the EPSRC [EP/L015234/1] via the LSGNT Centre for Doctoral Training.

\section{Toric background}

\subsection{Notation and assumptions}

We are interested in toric varieties constructed as GIT quotients of a vector space $V$ by a torus $T$. We specify the data of the torus action as a complex of lattices
\beq{eq.LtoN}
 0 \To L \stackrel{Q^\vee}{\To} \Z^n  \stackrel{A}{\To} N \To 0 \eeq
 or its dual:
\beq{eq.LtoNdual} 0 \To M \stackrel{A^\vee}{\To } \Z^n \stackrel{Q}{\To} L^\vee \To 0 \eeq
Here:
\begin{itemize}
\item $L$ is the lattice of 1-parameter subgroups of the torus $T$, so $T=L_{\C^*}$. 
\item $\Z^n$ is the lattice of Laurent monomials on $V$, \emph{i.e.}~$V=\mathrm{Spec}[\N^n]$ for the submonoid $\N^n\subset \Z^n$.%\footnote{Technically we could replace $\N^n$ with a free monoid $\cM$ and $\Z^n$ with $\cM_\Z$ but then $\cM$ has a canonical basis up to permutations so little is gained.}
\item $Q$ is the \emph{weight map}. The images $q_i= Q(e_i)$ of the standard basis vectors are the \emph{weights} of the action. 
\item $N$ is the cokernel of $Q^\vee$ modulo torsion.
\item $M$ is the kernel of $Q$ and the dual of $N$.
\item $A$ is the \emph{ray map}. The images $a_i=A(e_i)$ are the \emph{rays}. 
\end{itemize}

By definition $A$ is surjective and $A^\vee$ is injective. We will always assume that $Q^\vee$ is injective, so $Q$ is surjective modulo torsion - this is the assumption that generic points of our GIT quotient stacks do not have infinite isotropy groups.  It follows that \eqref{eq.LtoN} and \eqref{eq.LtoNdual} are exact apart from a possible torsion group $L^\vee/\Im Q \cong \Ker A / \Im Q^\vee$.

% We'll sometimes WHEN?? make the additional assumption that
%\beq{eq.nogerbe}Q \mbox{ is surjective} \eeq
%\emph{i.e.} $N$ is exactly the cokernel of $Q^\vee$, \emph{i.e.}  \eqref{eq.LtoN} and \eqref{eq.LtoNdual} are exact.  This says that $T\to GL(V)$ is injective so generic points in the quotient stacks have no stabilizer. Note that without this assumption knowing $A$ is not enough to recover $Q$.
%*********** I HAVEN'T FOUND ANYWHERE WE NEED THIS ASSUMPTION

\pgap

A \emph{stability condition} is an element of $L^\vee_\R$. A choice of stability condition $\theta$ defines a semi-stable locus in $V$ and hence a GIT quotient, which for us means the quotient stack:
$$ X_\theta = [ V^{ss}_{\theta} \, / \, T ] $$
We'll generally only be interested in quotients with respect to generic $\theta$, in which case $X_\theta$ is at worst a DM stack.  We'll also refer to these generic GIT quotients as the \emph{phases} of the GIT problem. Each phase is a toric orbifold and has a corresponding fan in $N$. The rays of this fan are always (some subset of) the $a_i$'s, hence the name. The higher dimensional cones change depending on the phase.

\begin{rem}\label{rem.Qnotsurj} If the weight map $Q$ has some finite cokernel then the representation $T\to GL(V)$ has a finite kernel, so the GIT quotients $X_\theta$ have finite isotropy groups at all points. We need to allow this possibility, since even if it doesn't apply to our initial toric variety $X$ it can happen for the smaller-dimensional varieties $Z$ that appear in wall-crossing.

 Note that in this situation $A$ does not determine $Q$. There is a theory of \emph{stacky fans} but which solves this issue but we won't need it because for us $Q$ is the fundamental piece of data.
\end{rem}
 
The space of stability conditions has a wall-and-chamber structure whose chambers correspond to phases.  If we consider all (non-empty) GIT quotients we get a fan in $L^\vee$ called the \emph{secondary fan} - the top-dimensional cones correspond to non-empty phases and the lower-dimensional cones correspond to non-generic GIT quotients. The rays of the secondary fan include those generated by the weights $q_i$, but in general there more rays than this.  Corresponding to the secondary fan is a toric variety, the \emph{secondary toric variety} $\mfF$.

% In fact $\mfF$ is more naturally an orbifold; there is a way to promote the secondary fan to a stacky fan, and then $\mfF$ is a toric stack \cite{Alex Sec 3.4 or C. Diemer, L. Katzarkov, G. Kerr.}. But this is not relevant for the current paper.

\subsection{The Calabi-Yau case}\label{sec.CYcase}

An important special case is when the torus $T$ acts through $SL(V)$, which implies that each phase is Calabi-Yau.

In terms of the toric data, the Calabi-Yau case is when the sum of the weights $q_i$ is zero. Equivalently, the rays $a_i$ are all contained in (and hence affinely span) an affine hyperplane of height 1. In this case is helpful to consider the polytope
$$\Pi \subset N_\R $$
given by the convex hull of the rays. Each phase corresponds to a fan in $N$, which when intersected with the affine hyperplane determines a decomposition of $\Pi$.  These decompositions are exactly the \emph{coherent triangulations}, \emph{i.e.}~triangulations induced by a piece-wise linear function. 

\subsection{Higgs and Coloumb GIT problems}\label{sec.higgsandcol}

From our original GIT problem $T\actson V$ we will often extract a smaller GIT problem involving some subset of the toric data, either by picking a subset of the weights, or a subset of the rays. The two main ways this will happen are:
\begin{enumerate}
\item Suppose $W\subset L^\vee_\R$ is a wall in the secondary fan, normal to some 1-parameter subgroup $\lambda \in  L$. Then we can consider the subset of weights which are orthogonal to $\lambda$, \emph{i.e.} which lie in the subspace $\langle W\rangle$.  

\item In the Calabi-Yau case we can choose a face $\Gamma\subset \Pi$ of the toric polytope, and consider the set of rays lying in this face.
\end{enumerate}

Formally, suppose we pick a subset $\cS\subset \{1, ..., n\}$. We can view $\cS$ as a subset of the standard basis vectors $\{e_1, .., e_n\}$ in $\Z^n$ so there is a corresponding set of rays $A(\cS) \subset N$. We set $N_\cS\subset N$ to be the sublattice spanned by $A(\cS)$, write $A_\cS: \Z^\cS \to N_\cS$ for the restriction of $A$, and set $L_\cS= \Ker A_\cS$. Then we get a GIT problem:
$$ L_\cS \stackrel{Q_\cS^\vee}{\To} \Z^\cS  \stackrel{A_\cS}{\To} N_\cS  $$
We'll refer to this as the \emph{Coloumb GIT problem} associated to the subset $\cS$.
\pgap

Alternatively we pick a subset $\cT \subset \{1, ..., n\}$ and consider the corresponding set of weights $Q(\cT) \subset L^\vee$. We define $L^\vee_\cT$ as the primitive sublattice generated by these weights
$$L_\cT^\vee = L^\vee  \cap \langle Q(\cT)\rangle_\R \;\subset L^\vee_\R $$
and we get a GIT problem:
$$ M_\cT \stackrel{A_\cT^\vee}{\To} \Z^\cT \stackrel{Q_\cT}{\To} L_\cT^\vee $$
We'll call this the \emph{Higgs GIT problem} associated to $\cT$.  Note that $Q_\cS$ is by definition surjective but $Q_\cT$ might not be (\emph{c.f.}~Remark \ref{rem.Qnotsurj}).
\pgap

Our `Higgs' and `Coloumb' terminology is based on the `Higgs GLSM' and `Coloumb GLSM' from \cite{APW}, which are related to the Higgs and Coloumb branches of the vacuum moduli space at singular values of the FI parameters.

\begin{rem}\label{rem.CCY} If our original GIT problem is Calabi-Yau then the Coloumb GIT problem is also Calabi-Yau for any subset $\cS$. But the Higgs GIT problems may not be. 
\end{rem}

\section{Semi-orthogonal decompositions for toric varieties}\label{sec.SODs}

\subsection{Crossing a single wall}\label{sec.singlewall}

Fix a toric GIT problem $T \actson V$. Let $C_+$ and $C_-$ be two adjacent chambers of the secondary fan separated by a wall $W$, and labelled such that $C_+$ lies on the same side of $W$ as the character $\det(V)$. Let $\lambda_W\in L$ be the primitive 1-parameter subgroup normal to this wall, oriented such that 
 $$ \kappa = (\det V)(\lambda_W) \geq 0 $$
\emph{i.e.}~$C_+$ lies on the $\lambda_W>0$ side.  Write $X_\pm$ for the phases corresponding to these two chambers. 

For this wall we have a Higgs GIT problem as described in Section \ref{sec.higgsandcol}. Let $\cT$ be the indexing set for the weights orthogonal to $\lambda_W$, so $Q(\cT)$ are all the weights lying in the subspace $\langle W \rangle$. The vector space corresponding to $\Z^\cT$ is the fixed subspace $V^{\lambda_W} \subset V$. Also $Q(\cT)$ necessarily span $\langle W \rangle$, so $L_\cT^\vee$ is exactly the orthogonal to $\lambda_W$, \emph{i.e.} it's the character lattice of $T/\lambda_W$. Hence this Higgs GIT problem is just describing the action of $T/\lambda_W$ on $V^{\lambda_W}$.

The secondary fan for this Higgs GIT problem lives in the vector space $\langle W \rangle$ and the cone $W$ lies in some chamber of it. We write $Z$ for the corresponding phase. 

\begin{thm}\cite[Theorem 5.2.1]{BFK}\label{thm.toricSOD}
We have a semi-orthogonal decomposition
$$D^b(X_+) = \SOD{ D^b(X_-), \; D^b(Z), \;...,\; D^b(Z) } $$
where $\kappa$ copies of $D^b(Z)$ occur.
\end{thm}

\begin{rem} This theorem is an application of the general theory of `windows' relating GIT and derived categories \cite{BFK, HL, Seg}, which applies to a general GIT quotient of a variety by a reductive group. However, in the current state-of-the-art you cannot use this theory to compare two \emph{different} GIT quotients unless you assume that the wall-crossing is of a particularly simple form. Which these ones are.
\end{rem}

\begin{rem} If $\det(V)$ lies on the wall then $\kappa=0$ and the theorem states that $D^b(X_+)$ and $D^b(X_-)$ are equivalent. This is a toric flop.
\end{rem}

\begin{eg}\label{eg.Pn} If we consider the standard action of $\C^*$ on $\C^{n+1}$ then $X_-=\varnothing$ and we get
$$D^b(\P^n) = \SOD{D^b(pt), \; ..., \; D^b(pt)} $$
which recovers Beilinson's result that $\P^n$ has full exceptional collection of length ${n+1}$. 
\end{eg}

\begin{rem}\label{rem.Orlov} If $X_+$ happens to be a blow-up of $X_-$ then Theorem \ref{thm.toricSOD} recovers Orlov's blow-up formula for this toric situation. It's possible to formulate the theorem more generally in such a way that it directly generalizes Orlov's result.
\end{rem}

\subsection{The algorithm}

Theorem \ref{thm.toricSOD} immediately suggests the following recursive algorithm for decomposing the derived category of a phase $X$:
\begin{enumerate}
\item Starting at the chamber for $X$ we cross through a sequence of walls, always moving away from $\det(V)$. At each wall we refine our decomposition. 
\item We stop when we reach a \emph{minimal phase} where no further such wall-crossings are possible.
\item Every factor occuring in this decomposition is the derived category of a phase of a smaller GIT problem, so we can apply this algorithm to each factor. 
\end{enumerate}

Note that a phase is minimal if $-(\det V)$ lies in the closure of that chamber, or equivalently if the canonical bundle of that phase is nef. 

\begin{rem} If $X$ is projective then you can use this algorithm to recover Kawamata's result \cite{Kaw} that a projective toric variety has a full exceptional collection \cite[Thm 5.2.3]{BFK}. This is because the minimal phase will be empty (as in Example \ref{eg.Pn}),  and moreover the minimal phase is empty in every Higgs GIT problem that occurs in the algorithm.
\end{rem}

In this paper we are more interested in quasi-projective examples. 

\begin{eg}\label{eg.algorithm}
Take $V=\C^6$ and quotient by $(\C^*)^2$ using the following matrix of weights:
$$\begin{pmatrix}
 1 & 1 & -1 & 0 & 0 & 0 \\ 
0 & 0 & 1 &1 & 1 & -1 \end{pmatrix} $$
Observe that $\det(V) = (1,2)^\top$. This GIT problem has four phases and the secondary fan is drawn in Figure \ref{fig.phases}. The phases are:
\begin{enumerate}
\item[(1)] $X_1=\A^4$. This is the unique minimal phase.
\item[(2)] $X_2=\cO(-1)_{\P^3}$, the total space of the tautological line bundle on $\P^3$.
\item[(3)] $X_3=\cO(-1)_{\P^1}\times \A^2$.
\item[(4)] $X_4 = \cO(-1)_P$, the total space of the relative $\cO(-1)$ line bundle over the projective bundle $P=\P(\cO^{\oplus 2}\oplus \cO(-1)) \to \P^1$.
\end{enumerate}

\begin{figure}[!tbp]
\begin{tikzpicture}[scale =1.4]
\draw[thick] (0,0) --(2,0);
\draw[thick] (0,-1)--(0, 1.3);
\draw[thick] (0,0)--(-1.5,1.3);
\node at (-.6,-.2) {$(1)$};
\node at (-.9,-.7) {$\A^4$};
\node at (-.3,.8) {$(2)$};
\node at (-.7,1.4) {$\cO(-1)_{\P^3}$};
\node at (.5,-.4) {$(3)$};
\node at (1.2,-.9) {$\cO(-1)_{\P^1}\times \A^2$};
\node at (.5,.5) {$(4)$};
\node[right] at (.8,1) {$\cO(-1)_{P}$};
\end{tikzpicture}
\caption{}\label{fig.phases}
\end{figure}

Firstly we decompose $D^b(X_2)$ by crossing the wall into chamber (1). The 1-parameter subgroup for this wall is $(1,1)$ so $\kappa=3$. The Higgs GIT is $\C^*\actson \C$ with weight 1, and $Z$ is the non-empty phase $Z=pt$. Hence Theorem  \ref{thm.toricSOD} in this case says
$$D^b(X_2) = \SOD{D^b(X_1),\; D^b(pt),\; D^b(pt),\; D^b(pt) } $$
which is an instance of Orlov's blow-up formula (see Remark \ref{rem.Orlov}). 

To make the rest of this example more readable we'll write this SOD and all following ones in the compressed form:
$$D^b(X_2) = \SOD{X_1,\; pt,\; pt,\; pt } $$
For this phase no futher refinements are possible, and the algorithm is finished.
\pgap

Next we apply the algorithm to phase 4. Let us choose to cross to phase (2) and then to phase (1).  The wall-crossing between (2) and (4) is again a blow-up, it blows up the codimension 2 subvariety $\cO(-1)_{\P^1}$. So crossing both walls gives:
$$D^b(X_4) = \SOD{ X_2, \; \cO(-1)_{\P^1}} = \SOD{ X_1, \; pt, \; pt, \; pt,\; \cO(-1)_{\P^1}} $$
We are not yet finished, because we can still apply the algorithm to the factor $D^b(\cO(-1)_{\P^1})$. But this variety is just the blow-up of $\A^2$ at the origin, so the next refinement is:
\beq{eq.SOD1}D^b(X_4) =  \SOD{ X_1, \; pt, \; pt, \; pt,\; \A^2, \; pt} \eeq
No further refinements are possible. 
\pgap

What happens if we make a different choice? We could instead have crossed to phase (3) before crossing to phase (1). The crossing $(1)\leadsto (3)$ blows up a plane, and the crossing $(3)\leadsto(4)$ blows up a $\P^1$, so crossing these walls gives the decomposition:
$$D^b(X_4) = \SOD{ X_3, \; \P^1, \; \P^1 } = \SOD{ X_1, \; \A^2, \; \P^1, \; \P^1 } $$
The factor $D^b(\P^1)$ can be split into two exceptional objects (as in Example \ref{eg.Pn}) so the final step is:
\beq{eq.SOD2}D^b(X_4) =  \SOD{ X_1, \; \A^2, \; pt, \; pt, \; pt, \; pt} \eeq
Note that in this example the quotienting torus is $(\C^*)^2$, and for each phase we needed to apply the recursive algorithm (at most) two times. For rank $r$ it would need $r$ applications. 
\end{eg}

\pgap

In the preceding example  we noticed that when decomposing $D^b(X_4)$ we had two choices, since there were two possible paths from chamber (4) to chamber (1).\footnote{By `path' we really mean a sequence of adjacent chambers.} In the second step of the algorithm there was no such choice, since the Higgs GIT problems were all rank 1 and had only two chambers. In a higher rank example there will be many more choices because we  need to choose a path at every step except the last one.

 However, examining the decompositions \eqref{eq.SOD1} and \eqref{eq.SOD2} that resulted from our two paths we can see evidence of our J\"ordan-Holder property - the decompositions are different, but the multiplicities of the `irreducible factors' agree. To state this precisely we need to think about what these `irreducible factors' really are.

\subsection{Relevant subspaces}\label{sec.relevantsubspaces} Recall that our initial GIT problem is given by a weight matrix $Q: \Z^n \to L^\vee$ specifying a torus action $T\actson V$. At any step in the algorithm the Higgs GIT problem arises as the fixed subspace $V^{T'}$ for some sub-torus $T'\subset T$, with  a corresponding sublattice $L'\subset L$. The weights $q_1,..., q_h$ of this Higgs GIT problem are those weights which are orthogonal to $L'$, and they always span the subspace $(L')^\perp_\R \subset L^\vee_{\R}$. 

The `irreducible factors' of our decompositions are the derived categories of the minimal phases of each such Higgs GIT problem. However, some of these minimal phases will be empty.  Since the stability condition that produces the minimal phase is 
$$-\det(V^{T'}) = -\sum_{i=1}^h q_i\quad \in (L')^\perp $$
 we get a non-empty minimal phase iff  the vector $-\sum q_i$ lies in the cone spanned by $q_1,..., q_h$. 

%\begin{defn} We call a subspace $H\subset L^\vee_{\R}$ \emph{relevant} if the cone spanned by the weights lying in $H$ is the whole of $H$.
%\end{defn}

\begin{defnlem}\label{defn.relevant} Let $H\subset L^\vee_{\R}$ be a subspace, let $q_1,..., q_h$ be the weights lying in $H$, and let $\sigma_H\subset H$ be the cone spanned by these weights. We call $H$ \emph{relevant} if one of the following two equivalent conditions hold:
\begin{enumerate}
\item[(i)] $\sigma_H$ is the whole of $H$.
\item[(ii)] $H$ is is spanned by $q_1,..., q_h$ and also $-\sum q_i \in \sigma_H $.
\end{enumerate}
\end{defnlem}
\begin{proof} 
Obviously (i) implies (ii). Conversely (ii) implies that $-q_i\in  \sigma_H$ for all $i$, so if the $q_i$'s span $H$ then any vector in $H$ can be written as a positive linear combination of them. 
\end{proof}

 Clearly there can only be finitely-many relevant subspaces. We allow $H=0$ (which is always relevant) and $H=L^\vee_{\R}$ (which might not be). A 1-dimensional relevant subspace is a line which has weights on both its rays. 

The relevant subspaces index the `irreducible factors' in our semi-orthogonal decompositions.  Each one defines a Higgs GIT problem with a non-empty minimal phase $Z_H$, and the corresponding factor is $D^b(Z_H)$. 

\begin{eg} In Example \ref{eg.algorithm} there are three  relevant subspaces: the whole of $\R^2$, the vertical axis, and the origin. They contribute the factors $D^b(\A^4)$, $D^b(\A^2)$ and $D^b(pt)$ respectively. 
\end{eg}

\begin{thm}\label{thm.JH}
 Let $X$ be a phase of a toric GIT problem and let $H$ be a relevant subspace. The multiplicity of $D^b(Z_H)$ in the semi-orthogonal decomposition of $D^b(X)$ is independent of all choices of paths. 
\end{thm}

Presumably the different decompositions resulting from different choices of paths are always related by mutations, but we haven't checked this.

\begin{rem}\label{rem.mult} The actual value of the multiplicity of $D^b(Z_H)$ in $D^b(X)$ is determined algorithmically from the toric data. It would be interesting -  and probably helpful for Conjecture \ref{conj.mainconj} -  to find something like a closed-form expression for it. 

We don't know how to do this, except in the case when $H$ has codimension 1 when it follows easily from the discussion in Section \ref{sec.singlewall}. Let $\lambda_H$ be a primitive 1-parameter subgroup normal to $H$, oriented so that it pairs positively with the chamber for $X$, and set $\kappa = (\det V)(\lambda_H)$. Then the multiplicity of $D^b(Z_H)$ in $D^b(X)$ is
$$\max\{ \kappa, 0 \}$$
since the algorithm only tells us to cross $H$ if $\kappa>0$.
\end{rem}

\subsection{Proof of the main theorem}

We'll prove Theorem \ref{thm.JH} using the recursive structure of the algorithm to reduce to the rank 2 case, \emph{i.e.} when the GIT problem consists of $(\C^*)^2\actson V= \C^n$. In the rank 1 case the theorem is vacuous since there are no choices.

\begin{lem} \label{lem.rank2}Theorem \ref{thm.JH} holds in the rank 2 case. 
\end{lem}
\begin{proof}
If $\det(V)$ is the trivial character then all phases are derived equivalent and the theorem is vacuously true, so we can assume $\det(V)\neq 0$. For simplicity we assume that neither $\det(V)$ nor $-\det(V)$ lie on a wall, so there is a unique minimal phase 
and a unique `maximal' phase $X_{max}$, whose chamber contains $\det(V)$. In fact there could be up to two minimal or maximal phases, but crossing the walls between them is a derived equivalence and we can ignore it. If we start at any non-maximal phase then there are no choices to be made in the algorithm, but if we start at $X_{max}$ then we have exactly two  choices of paths to reach $X_{min}$. So the only thing to check is that these two choices produce the same multiplicities.

There are three classes of relevant subspace:

\begin{enumerate}[leftmargin=18pt, itemsep=5pt]

\item  $H=\C^2$. This is relevant iff $X_{min}$ is non-empty, in which case $D^b(X_{min})$ occurs in $D^b(X_{max})$ with multiplicity one for either choice of path. 

\item $H$ a line, both rays of which are walls. The Higgs GIT for $H$ has two non-empty phases, let  $Z_H$ be a minimal one and $Z_H'$ be the other one. 

By assumption $\det(V)$ doesn't lie on $H$, so if $\lambda_H$ is a primitive normal 1-parameter subgroup to $H$ then  $\kappa= |\lambda_H(\det V)|$ is strictly positive. The minimal and maximal chambers lie on opposite sides of $H$ so either choice of path crosses it; one choice contributes $\kappa$ copies of $D^b(Z_H)$ and the other contributes $\kappa$ copies of $D^b(Z_H')$. But the decomposition of $D^b(Z_H')$ includes exactly one copy of $D^b(Z_H)$ so either way the multiplicity of $D^b(Z_H)$ in $D^b(X_{max})$ is $\kappa$.

\item $H=\{0\}$. This contributes the factor $D^b(V^T)$, the subspace of $V$ fixed by the whole torus. 

Consider a line  $l \subset L_\R^\vee$ containing at least one weight, let $q_l$ be the sum of the weights on this line, and let $\mu_l = |q_l|$ be the lattice length of $q_l$. There are two possibilities:
\begin{enumerate} \item There are weights on both rays of $l$. Then $l$ is a relevant subspace as in case (2), both rays are walls and the Higgs GIT has a non-empty minimal phase $Z_l$. The derived category of the other phase $Z_l'$ decomposes into one copy of $D^b(Z_l)$ and $\mu_l$ copies of $D^b(V^T)$.
\item There are only weights on one ray so only that ray is a wall. The Higgs GIT has an empty phase and the other phase decomposes into $\mu_l$ copies of $D^b(V^T)$.
\end{enumerate}
In either case only one of our two paths will pick up any factors of $D^b(V^T)$ from this line $l$; it's the path that crosses $l$ on the same side as $q_l$, and the number of such factors it picks up is
$$ \mu_l \kappa_l = \mu_l |\lambda_l(\det(V))| $$
where $\lambda_l$ is a primitive 1-parameter subgroup normal to $l$. So we may as well assume that each such line contains only a single weight $q_i=q_l$, and hence only that ray of the line is a wall. 

Now fix an orientation on our lattice $L^\vee$. Since the lattice is rank 2 this is the same as a unit symplectic form $\omega$. This means that for the wall through $q_i$ we can produce a primitive normal subgroup by setting  $\lambda = \omega(\hat{q_i}, -)$ where $\hat{q_i}$ is a primitive vector in the direction of $q_i$. With this choice one of our paths always crosses walls in the direction of increasing $\lambda$ and the other path always crosses walls in the direction of decreasing $\lambda$. So if the first path crosses the rays through $q_1, ..., q_s$ and the second path crosses the rays through $q_{s+1}, ..., q_n$ then the equality we want to show is:
$$\sum_{i=1}^s \mu_i \lambda_i(\det(V)) = - \sum_{i=s+1}^n  \mu_i \lambda_i(\det(V)) $$
But this is true since
$$ \sum_{i=1}^n \mu_i\lambda_i = \sum_{i=1}^n \omega( q_i, - ) = \omega(\det(V), -) $$
and $\omega(\det(V), \det(V))=0$. 
\end{enumerate}
\end{proof}

Now suppose have have a higher rank problem, and we choose a phase $X$ corresponding to a chamber $C_X$. To run the algorithm we first pick a path from $C_X$ to the chamber for a minimal phase, $C_{min}$, always moving away from $\det(V)$.
It doesn't matter which minimal phase we pick since moving between them is a derived equivalence. But there might be many possible paths from $C_X$ to $C_{min}$. 

To visualize this clearly pass from the secondary fan to the dual `secondary polytope' in $L_\R$. This has a vertex for each chamber, an edge for each wall, and higher-dimensional faces for fans of higher codimension. The element $\det(V)\in L^\vee$ defines a linear function on the polytope and induces a direction on (most of) the edges, since we never allow this function to increase when we traverse an edge.

Choose a path,  \emph{i.e.} a directed sequence of edges, between the vertices $c_X$ and $c_{min}$ corresponding to the chambers $C_X$ and $C_{min}$. Now pick a polygon $P$ (a two-dimensional face) in the secondary polytope which meets our path; let's say the path meets $P$ at some vertex $c_1$, traverses some edges of the polgon, then leaves it again at $c_2$. If the remaining edges in $P$ also happen to form a directed path then we can produce a new path from $c_X$ to $c_{min}$ by choosing to go the other way around $P$. This is possible iff $c_1$ maximizes $\det(V)$ among vertices of $P$ and $c_2$ minimizes $\det(V)$.

 Let's call this kind of operation on paths a \emph{simple modification}.

\begin{lem}\label{lem.moves} Any two paths from $c_X$ to $c_{min}$ are connected by a sequence of simple modifications. 
\end{lem}
\begin{proof}
The secondary polytope is a cell decomposition of an $(r-1)$-sphere, and the subset where $\det(V)\leq \det(V)(c_X)$ is a decomposition of a disc. Any two paths in this disc from $c_X$  to $c_{min}$ are homotopic, and given such a homotopy we can move it orthogonally to $\det(V)$ until it lies in the 2-skeleton of the polytope.  We then have a collection of polygons whose boundary is the union of our two paths, with $\det(V)\leq \det(V)(c_X)$ everywhere. We just need to show that we can perform a simple modification to one of our paths at one of these polygons; then the result follows by induction. 

Choose one of our paths. Let $P_1$ and $P_2$ be the first two polygons this path meets and let $c_1$ be the vertex where the path switches between them. By considering the edge between $P_1$ and $P_2$, and remembering that $\det(V)$ is a linear function on these affine polygons, we can see that either $c_1$ minimizes $\det(V)$ in $P_1$, or else it maximizes $\det(V)$ in $P_2$. In the first case we can modify the path at $P_1$. In the second case we move to the next pair of polygons along the path and repeat the argument. Note that $c_{min}$ certainly minimizes $\det(V)$ in the final polygon so the algorithm terminates there if not before.
\end{proof}

\begin{proof}[Proof of Theorem \ref{thm.JH}]
Pick two paths from the chamber for $X$ to the chamber for a minimal phase, always moving away from $\det(V)$. By Lemma \ref{lem.moves} it's enough to deal with the case when our two paths are related by a simple modification; this means they agree except at a single codimension-two cone $\Lambda$ in the secondary fan where they travel opposite ways around. For every Higgs GIT problem that our paths encounter we also need to make choices, but those GIT problems have lower rank so by induction we can assume that those choices do not matter. 

 Let $U\subset V$ be the semi-stable locus for a character lying on our codimension-two cone. Then we have a GIT problem $T\actson U$ whose phases are exactly those phases of $T\actson V$ whose chambers are adjacent to the cone. This new GIT problem is `non-linear' in that $U$ is not a vector space, and there is an important sense in which it is rank two. If we let $L'\subset L$ be the rank 2 sublattice normal to our codimension-two cone, and $T'\subset T$ be the corresponding subtorus, then only subgroups lying in $T'$ can have fixed points in $U$. It follows that the GIT fan for $T\actson U$ is just the GIT fan for $T'\actson U$, pulled-back via the projection $L^\vee_\R \to (L')^\vee_\R$. 

So in the region where our two paths differ we can think of them as paths in the GIT fan for $T'\actson U$. And since they are different they both must start in a maximal chamber and end in a minimal chamber. 

Now consider the linear GIT problem $T'\actson V$. The GIT fan for $T'\actson U$ is a coarsening of the one for $T'\actson V$; every wall of the former is a wall of the latter, but not necessarily vice-versa since a subgroup $\lambda\subset T'$ could have fixed points in $V$ but none in $U$. However, from the point-of-view of our algorithm there is no harm in regarding every wall for $T'\actson V$ as corresponding to a wall for $T\actson U$ -  it just happens that some of them will be `fake walls' where the semi-stable locus does not change. In the semi-orthogonal decomposition crossing a fake wall adds some number of copies of the zero category  $D^b(U^\lambda /\!/ T) = D^b(\varnothing)$. Note that if both rays of a line are a wall for $T'\actson V$ then either both give genuine walls for $T\actson U$ or both are fake. Also if there are any fake walls then $U^{T'}$ is empty, which means that the codimension-two cone itself also contributes the zero category. 

If we include these zero categories then we have a bijection between the factors in the decomposition algorithms for $T'\actson V$ and for $T\actson U$, and their multiplicities agree since these depend only on the restriction of the character $\det(V)$ to the subtorus $T'$. Hence the result follows from Lemma \ref{lem.rank2}.

\end{proof}

\section{FI parameter spaces and discriminants}

In this section we consider a \emph{Calabi-Yau} GIT problem $T\actson V$ where $T$ acts through the subgroup $SL(V)$. This has  a different flavour to the previous section, since all phases are Calabi-Yau and every wall-crossing is a derived equivalence, so no semi-orthogonal decompositions occur. Instead (as discussed in the introduction) we focus on autoequivalences of the phases and relate these to the fundamental group of the FI parameter space.

\subsection{Spherical functors}\label{sec.spherical}

Let $T\actson V$ be a Calabi-Yau toric GIT problem. Let $X_+$ and $X_-$ be two phases coming from two adjacent chambers $C_+$ and $C_-$, separated by a wall $W$. Let $Z$ be the phase of the associated Higgs GIT problem for a character lying on $W$. 

Since $\det(V)=0$, Theorem \ref{thm.toricSOD} tells us that $D^b(X_+)$ and $D^b(X_-)$ are equivalent. However, what the theory actually gives us is a countable set of equivalences
$$\Phi_i : D^b(X_+) \isoto D^b(X_-)$$
indexed by the integers. They are related by the Picard groups of $X_+$ and $X_-$. 

\begin{thm}\cite[Prop. 3.4]{HLSh} 
 There is a spherical functor
$$F: D^b(Z) \to D^b(X_+)$$
such that $\Phi_{1}^{-1}\Phi_0$ is the twist around $F$. 
\end{thm}

Recall that the \emph{twist} around $F$ is the endofunctor of $D^b(X_+)$ defined by the cone on the counit
$$T_F = [FR \to \id] $$
where $R$ is the right adjoint to $F$, and that the key property of a \emph{spherical} functor is that $T_F$ is an autoequivalence. See \cite{AL} for more detail on spherical functors. Note that this cone of functors makes sense since we can interpret it as a cone of Fourier-Mukai kernels (or insert the prefix `dg' where needed). 
\pgap

The variety $Z$ is toric - it's a phase of the Higgs GIT problem - but it will not usually be Calabi-Yau. So using the algorithm of Section \ref{sec.SODs} we can produce a semi-orthogonal decomposition:
\beq{eq.ZSOD}D^b(Z) = \SOD{ \cC_1,..., \cC_r} \eeq
Halpern-Leistner and Shipman observed that this implies:
\begin{enumerate}\item The restriction of $F$ to each piece gives a spherical functor $F_i: \cC_i \to D^b(X_+)$.
\item The twist $T_F$ factors as:
\beq{eq.twistfactors}T_F = T_{F_1}\circ... \circ T_{F_r} \eeq
\end{enumerate}
The formal result is \cite[Theorem 4.14]{HLSh} and it applies in this situation since the cotwist around $F$ is (up to a shift) the Serre functor on $D^b(Z)$. 
\pgap

The factors in the semi-orthogonal decomposition \eqref{eq.ZSOD} are indexed by the relevant subspaces in the Higgs GIT problem for $W$, but these are simply the relevant subspaces $H\subset L^\vee_\R$ which are contained in the hyperplane $\langle W \rangle$. 

\subsection{Discriminants}\label{sec.discriminants} We now recall some of the theory of discriminant loci developed by Gelfand--Kapranov--Zelevinsky \cite{GKZ}. 

Recall that our GIT problem is specified by a sequence of lattices, exact modulo torsion, or its dual:
$$L \stackrel{Q^\vee}{\To} \Z^n \stackrel{A}{\To} N$$
$$M \stackrel{A^\vee}{\To} \Z^n \stackrel{Q}{\To} L^\vee$$
From this point on we need to make two mild additional assumptions:
\begin{enumerate}\item We assume that the rays $A(e_i)$ are all distinct. We need this because for \cite{GKZ} $A$ is a subset of $N$. This excludes 1-parameter subgroups acting with weights $(0,...,0, 1,-1,0,...,0)$ but these are very uninteresting from a wall-crossing perspective.
\item We assume the weights $Q(e_i)$ are all non-zero. This is just for simplicity. A zero weight just contributes a factor of $\A^1$ to each phase.
\end{enumerate}

 Tensoring our lattices by $\C^*$ gives two exact sequences of tori:
$$L_{\C^*} \stackrel{Q^\vee_\times}{\To}(\C^*)^n \stackrel{A_\times}{\To} N_{\C^*}$$
$$M_{\C^*} \stackrel{A^\vee_\times}{\To} (\C^*)^n\stackrel{Q_\times}{\To} L^\vee_{\C^*}$$

The map $A^\vee_\times$ provides us with $n$ characters of the torus $M_{\C^*}$. If we pick a vector of coefficients $a\in \C^n$ we can take a linear combination of these characters, this gives us a Laurent monomial:
\al{W_a : M_{\C^*} &\to \C \\
x&\mapsto \langle a, A^\vee_\times(x) \rangle}
In explicit coordinates this means
$$W_a  =\sum_{i=1}^n a_i \prod_{t=1}^m X_t^{A_{it}} $$
where $X_1,..., X_m$ are coordinates on $M_{\C^*}$. This is the \emph{Hori-Vafa mirror} to our toric GIT problem (or abelian GLSM); it's a family of Landau-Ginzburg models parametrized by $a$.

Since our GIT problem is Calabi-Yau we can choose co-ordinates such that the first column of $A$ is entirely 1's, hence 
$$W_a = X_1 \widetilde{W}_a $$
where $X_1$ doesn't appear in $\widetilde{W}_a$. 

For a generic $a$ the zero locus $W_a$ will be a smooth hypersurface in $M_{\C^*}$. Consider the subset of non-generic $a$, \emph{i.e.}
$$ D_A = \left\{ a\in \C^n, \; \exists x\in M_{\C^*} \mbox{ such that } W_a(x)=0\mbox{ and } dW_a(x)=0 \right\}$$
This, or perhaps its closure, is the \emph{discriminant locus} of the family $W_a$.  This definition is the correct one for general $A$; since we're in the Calabi-Yau case the first condition is redundant as $\partial_{X_1}W_a=0$ implies $W_a=0$. 

The closure of $D_A$ is an affine variety, which is always irreducible and usually a hypersurface \cite[Ch.~9]{GKZ}. To understand why this is true observe that $D_A$ is a cone so there is an associated projective variety in $\P^{n-1}$. It's not hard to compute that its projective dual is the closure of the image of $M_{\C^*}$ in $\P^{n-1}$, which is evidently irreducible. But the projective dual to an irreducible variety is always irreducible, and usually a hypersurface \cite[Ch.~1]{GKZ}.  If $D_A$ is a hypersurface then we denote its defining polynomial by $\Delta_A$. 

As well as being a cone $D_A$ is invariant under rescaling the $X_i$ variables, \emph{i.e.} it is invariant under the action of the torus $M_{\C*}$ on $\C^n$. We can replace $D_A$ with the open subset $D_A\cap (\C^*)^n$ - if $D_A$ is a hypersurface this loses no information - and then the quotient by $M_{\C^*}$ is a subvariety:
$$\nabla_A \subset L^\vee_{\C*}$$
$D_A$ is a hypersurface iff $\nabla_A$ is, and in this case $\Delta_A$ is really a function on $L^\vee_{\C^*}$.

%One of the main results of \cite[Ch. 10, Theorem 1.4]{GKZ} is that (assuming $\nabla_A$ is a hypersurface) the Newton polytope of $\Delta_A$ is the toric polytope of$\mathfrak{F}$. This means that $\overline{\nabla}_A$ avoids the toric fixed points, and the tropicalization of $\nabla_A$ is the secondary fan.  In this sense $\overline{\nabla}_A$ is the natural compactification of $\nabla_A$. 

%Calabi-Yau condition implies that $0$ can be the only critical value so $W_a(x)$ is redundant. Projective-dually it implies that the image of $M_{\C^*}$ in $\C^n$ is a cone.

\subsubsection{Horn uniformization}\label{sec.Horn}

In the Calabi-Yau case there is a useful dominant rational map
$$\P L_\C  \dashrightarrow \nabla_A$$
 called Horn uniformization, given by: 
$$[\lambda] \mapsto Q_\times \circ Q^\vee(\lambda) $$
In explicit co-ordinates this says:
$$\lambda_1\!:\! ... \!:\! \lambda_r \; \mapsto \; \left( \prod_{i=1}^n \Big( \sum_{k=1}^r Q_{ik}\lambda_k\Big)^{Q_{i1}}\!\!, \; \; ...\; \; , \; \;   \prod_{i=1}^n \Big( \sum_{k=1}^r Q_{ik}\lambda_k\Big)^{Q_{ir}}      \right) $$

\begin{eg}\label{eg.rank1discriminant} Suppose $L=\Z$ has rank one, and write $(q_1, ..., q_n)$ for the vector of weights. Then by the above $\nabla_A$ consists of the single point
$$ q_1^{q_1}... q_n^{q_n} \in \C^*$$
(recall we are assuming that no weights are zero). In particular $\nabla_A$ is a hypersurface and non-empty.
\end{eg}

Let's explain why this works. We have:
$$\partial_{X_s} W_a =\frac{1}{X_s}\sum_{i=1}^n a_i A_{is}\prod_{t=1}^m X_t^{A_{it}} $$
Invariantly, for a fixed $x\in M_{\C^*}$ this says that $dW_a(x)$ is the linear map
$$ dW_a(x) : M_\C  \To \C $$
given by composing:
$$M_\C \stackrel{x^{-1}}{\To} M_\C \stackrel{A^\vee}{\To} \C^n \stackrel{A^\vee_\times(x)}{\To} \C^n \stackrel{a}{\To}\C$$
Here the first map is the action of the element $x^{-1}\in M_{\C^*}$ on $M_\C$, and similarly for the third map. So  $dW_a$ has a critical point at $x$ iff $a\circ A^\vee_\times(x)$ annihilates $M_\C$, \emph{i.e.} iff 
$$a\circ A^\vee_\times(x) = Q^\vee (\lambda)$$
for some $\lambda \in L_\C$. So the image of the map
\al{ M_{\C^*} \times L_\C &\To \C^n\\
(x, \lambda) & \mapsto a= \big(A^\vee_\times(x)\big)^{-1}Q^\vee(\lambda) }
is the subset where $W_a$ has a critical point, and in the Calabi-Yau case this  is exactly $D_A$.  

Next we compose this with the quotient map $Q_\times: D_A \dashrightarrow \nabla_A$ and observe that
$$Q_\times(a) = \big(Q_\times A^\vee_\times(x)\big)^{-1} Q_\times Q^\vee(\lambda) = Q_\times Q^\vee(\lambda)$$
is independent of $x$, since $Q_\times A^\vee_\times(x) =1$. Hence $Q_\times \circ Q^\vee$ is a dominant rational map from $L_\C$ to $\nabla_A$. Finally, the Calabi-Yau condition implies that this map descends to $\P L_\C$. 
\pgap

If $\nabla_A$ is a hypersurface it has the same dimension as $\P L_\C$, and in this case Horn uniformization is a birational equivalence \cite[Ch.~9, Thm 3.3]{GKZ}. The inverse is the logarithmic Gauss map. 

\subsubsection{Components of the discriminant}

Recall that the convex hull of the rays $A(e_i)$ is a polytope $\Pi \subset N_\R$, which lies in an affine hyperplane of height 1. 

Choose a face $\Gamma$ of $\Pi$. Associated to this face there is a Coulomb GIT problem, as described in Section \ref{sec.higgsandcol}. We consider all the rays that lie in this face, and (abusing notation) write $\Gamma\subset\{1,..., n\}$ for the subset that indexes these rays. Then the Coulomb GIT problem is specified by an exact sequence of lattices
\beq{eq.LGamma}L_\Gamma \stackrel{Q^\vee}{\To} \Z^\Gamma \stackrel{A_\Gamma}{\To} N_\Gamma\eeq
where $N_\Gamma$ is the sublattice spanned by the face. 

We can define a discriminant locus associated to this face in the same way as we did for the whole polytope. For any vector of coefficients $a' \in \C^\Gamma$ there is a Laurent monomial $W'_{a'}$ on the torus $M^\Gamma_{\C^*}$, where $M^\Gamma$ is the dual lattice to $N_\Gamma$. To obtain $W'_{a'}$ from $W_a$ you just delete all the terms that don't correspond to rays on $\Gamma$, then since only some variables remain this function descends from $M_{\C^*}$ to the quotient $M^\Gamma_{\C^*}$.  Proceeding as before, we obtain a discriminant subset $D_\Gamma\subset \C^\Gamma$, a subvariety
$$ \nabla'_\Gamma \subset  (L^\vee_\Gamma)_{\C^*} $$
and its preimage:
$$ \nabla_\Gamma \subset  (L^\vee)_{\C^*} $$

\begin{rem}\label{rem.notaface} What we've just done works for any subset of the rays, not just the subsets corresponding to faces of $\Pi$. But the faces are the most important. Also note that the Coulomb GIT problems are all Calabi-Yau (Remark \ref{rem.CCY}) so we still have Horn uniformization.
\end{rem}

Roughly, we are interested in the union of these subvarieties over all faces of $\Pi$. However, some faces don't contribute anything. For example if $\Gamma$ is a simplex then $L_\Gamma=0$ so $\nabla_\Gamma$ must be empty; indeed it's easy to see that $D_\Gamma$ is just the origin in this case. 

More generally suppose $\Gamma$ contains a ray $A(e_i)$ which is linearly independent of the other rays in $\Gamma$. Then $D_\Gamma$ will be contained in the hyperplane $a'_i=0$ and hence $\nabla_\Gamma$ is empty. If we want to access $D_\Gamma$ then we should try deleting this ray $A(e_i)$; this will give us a subface $\Sigma \subset \Gamma$ with one less ray but with $L_\Sigma = L_\Gamma$. Then $D_\Sigma= D_\Gamma$ under the inclusion $\C^\Sigma \into \C^\Gamma$, but $\nabla'_\Sigma$ might be a non-empty subvariety of the torus $(L^\vee_\Sigma)_{\C^*} = (L^\vee_{\Gamma})_{\C^*}$. This observation leads to us to the following:

\begin{defn} \label{def.minimal} A subset $\cS\subset\{1,..., n\}$ is \emph{minimal} if, for all $i\in \cS$, the ray $A(e_i)$ is linearly dependent on the remaining rays $\{A(e_j), j \in \cS \setminus i\}$.

A face $\Gamma\subset \Pi$ is \emph{minimal} if the set of all rays lying in $\Gamma$ is indexed by a minimal subset.
\end{defn}
So a face $\Gamma$ is minimal iff we can remove any ray from $\Gamma$ without making the linear span smaller. Then we define:

\begin{defn} The \emph{discriminant locus} $\nabla\subset (L^\vee)_{\C^*}$ is the union of the subvarieties $\nabla_\Gamma$, for each minimal face $\Gamma\subset \Pi$ such that $\nabla_\Gamma$ is a hypersurface.
\end{defn}

The whole polytope $\Pi$ is minimal since we're assuming that there are no zero weights. If $\nabla_\Pi=\nabla_A$ is a hypersurface then we call it the \emph{principal component} of $\nabla$.

\begin{rem} This definition comes from \cite{GKZ}. It is not entirely clear to us why one disregards the subvarieties $\nabla_\Gamma$ which are not hypersurfaces. In the examples we've calculated it makes no difference, \emph{i.e.} each discriminant subvariety of higher codimension is contained in one which is a hypersurface. But we don't know if this is always true.
\end{rem}

If $\nabla_\Gamma$ is a hypersurface we write $\Delta_\Gamma$ for its defining polynomial, then the product of these cuts out the hypersurface $\nabla$.  Gelfand--Kapranov--Zelevinzky modify this by introducing some multiplicities $\mu_\Gamma$   and then taking the product
$$E_A = \prod_\Gamma (\Delta_\Gamma)^{\mu_\Gamma}$$
which they call the \emph{principal $A$-determinant} \cite[Ch.~10, 1.B]{GKZ}. The $\mu_\Gamma$'s are not relevant for us but there are two important theorems that they prove that are stated in terms of $E_A$. 

\begin{thm}\cite[Ch.~10, Thm 1.4]{GKZ}\label{thm.Newtonpolytope} The Newton polytope of $E_A$ is dual to the secondary fan. 
\end{thm}

In fact they give a more precise definition of the \emph{secondary polytope} $\check{\Pi}$ - which is in particular dual to the secondary fan - and their theorem is that the Newton polytope of $E_A$ is $\check{\Pi}$. There is a potential sign confusion here: the theorem is that the cones of the secondary fan are the same as the cones spanned by the inward normal vectors at each vertex of $\check{\Pi}$.

Recall that the secondary fan is the fan of the secondary toric variety $\mathfrak{F}$. This is a compactification of $L^\vee_{\C^*}$ so we can consider the closure:
$$\overline{\nabla}  \subset \mathfrak{F} $$
The theorem above suggests that this is a natural choice of compactification for $\nabla$. In particular it implies:

\begin{cor}\label{cor.missesfp} $\overline{\nabla}$ avoids all the toric fixed points in $\mathfrak{F}$.\end{cor}
\begin{proof} A fixed point is the origin in one of the toric charts. Each chart corresponds to a vertex of the Newton polytope of $E_A$, and when we write $E_A$ in that chart we get a non-zero constant term.
\end{proof}
\pgap

Recall also that phases of our GIT problem correspond to coherent triangulations of the polytope $\Pi$, meaning triangulations induced by a piece-wise linear function \cite{GKZ}. More generally a non-generic stability condition induces a \emph{coherent subdivision} of $\Pi$ where not all the pieces are simplices. Such a stability condition corresponds to a face of the secondary polytope $\check{\Pi}$ whose vertices are the phases refining this subdivision to a triangulation. 

Suppose we fix a coherent subdivision of $\Pi$, corresponding to a face $\check{\Gamma}\subset \check{\Pi}$. Now choose one of the pieces of the subdivision, it is some polytope $\Sigma_i \subset \Pi$.  As usual we abuse notation and also write $\Sigma_i\subset\{1,..., n\}$ for the indexing set of the rays appearing in this polytope. Associated to this subset $\Sigma_i$ we have a Coloumb GIT problem and a corresponding discriminant locus $\nabla_{\Sigma_i}\subset L^\vee_{\C^*}$ (see Remark \ref{rem.notaface}).  If $\Sigma_i$ is a simplex this discriminant locus is empty, so it's only worth considering the non-simplicial pieces of our subdivision.

Going further we can consider the principal determinant $E_{\Sigma_i}$, which we may view as a function on $L^\vee_{\C^*}$ by pulling-back under the projection $L^\vee\to L_{\Sigma_i}^\vee$. The zero locus of $E_{\Sigma_i}$ consists of the discriminant locus associated to $\Sigma_i$  as well as the discriminant loci coming from all the faces of $\Sigma_i$.

On the other hand, Theorem \ref{thm.Newtonpolytope} tells us that the face $\check{\Gamma}$ corresponds to some subset of the monomials appearing in $E_A$. Let us write $(E_A)_{\check{\Gamma}}$ for the sum of this set of monomials.

 \begin{thm}\cite[Ch.~10, Thm 1.12]{GKZ}\label{thm.EAmultiplies} For some positive integer multiplicities $\mu_i$ and some non-zero constant $\nu$ we have
$$(E_A)_{\check{\Gamma}} = \nu \prod_{i} \left(E_{\Sigma_i}\right)^{\mu_i}$$
where the product runs over the non-simplicial pieces of the subdivision. 
\end{thm}

In fact we only care about one special case of this theorem: the case when $\check{\Gamma}$ is an edge of $\check{\Pi}$. Such an edge connects two phases, and corresponds to a wall $W$ in the secondary fan. In the secondary toric variety $\mfF$ the phases correspond to toric fixed points, and the wall $W$ (or edge $\check{\Gamma}$) corresponds to a toric rational curve
$$C_W \subset \mfF$$ connecting the two fixed points. We discussed this in Section \ref{sec.multsintro}.

\begin{cor}\label{cor.onepointinCW} The discriminant locus $\overline{\nabla}$ intersects $C_W$ in exactly one point. 
\end{cor}
\begin{proof}
The intersection of $\overline{\nabla}$ with $C_W$ is the zero locus of the restriction $E_A|_{C_W}$ and this restriction is the sum $(E_A)_{\check{\Gamma}}$ of the monomials appearing in the edge $\check{\Gamma}$. This edge corresponds to a coherent subdivision of $\Pi$ which has exactly one non-simplicial piece $\Sigma$, having two possible triangulations.  By Theorem \ref{thm.EAmultiplies} the zero locus of $(E_A)_{\check{\Gamma}}$ agrees with the zero locus of $E_{\Sigma}$.

But the zero locus of $E_\Sigma$ is the discriminant locus for Coloumb GIT problem associated to $\Sigma$. This GIT problem has $\rank L_\Sigma = 1$ so by Example \ref{eg.rank1discriminant} its discriminant locus is a single point.
\end{proof}

In Lemma \ref{lem.intersection} below we will refine this result by identifying which components of $\overline{\nabla}$ can intersect with $C_W$.

\subsection{Faces and subspaces}

In Section \ref{sec.relevantsubspaces} we discussed \emph{relevant subspaces} in $L^\vee_\R$, these index the factors appearing in our SODs. In this section we show that relevant subspaces biject with minimal faces of the polytope $\Pi$; this is an elementary observation but crucial for formulating our conjecture.
\pgap

Consider a subset $\cS\subset \{1,..., n\}$ and its complement $\cS^c$. Let's consider the Coloumb GIT problem associated to $\cS$ and the Higgs GIT problem associated to $\cS^c$ (Section \ref{sec.higgsandcol}). These are related by the following diagram:
\beq{eq.higgscoloumb}
\begin{tikzcd}
M_{\cS^c} \arrow[r, "A_{\cS^c}^\vee"] \arrow[d] & \Z^{\cS^c} \arrow[r, "Q_{\cS^c}"]\arrow[d] & L_{\cS^c}^\vee \arrow[d]\\
M \arrow[r, "A^\vee"]\arrow[d] & \Z^n \arrow[r, "Q"] \arrow[d] & L^\vee \arrow[d] \\
M_\cS \arrow[r, "A_\cS^\vee"] & \Z^\cS \arrow[r, "Q_\cS"] & L_\cS^\vee
\end{tikzcd}
\eeq

The middle column is obviously exact, the other columns are exact modulo torsion. Let us also write 
$$H_{\cS^c} \;\subset L^\vee_\R $$
for the subspace spanned by $L_{\cS^c}^\vee$. 

As a special case we could consider a face of the polytope $\Pi$ and let $\Gamma$ be the indexing set for the rays on that face. Then we get an associated subspace $H_{\Gamma^c}\subset L^\vee_\R$. 

\begin{prop}\label{prop.minimaltorelevant} The map $\Gamma \mapsto H_{\Gamma^c}$ is a bijection between the minimal faces of $\Pi$ and the relevant subspaces  of $L^\vee_\R$. 
\end{prop}
\begin{proof}
 Recall from Definition \ref{def.minimal} that a subset $\cS\subset\{1,..., n\}$ is called minimal if the set of rays $A(\cS)\subset N$ has the property that every ray in $A(\cS)$ is linearly dependent on the remaining rays in $A(\cS)$. This is the statement that no basis vectors map to zero under the map $\Z^\cS \to L_\cS^\vee$, or equivalently that the only weights lying in $H_{\cS^c}$ are $Q(\cS^c)$.  

Conversely, pick a subspace $H\subset L^\vee_\R$ which is spanned by the weights it contains, and let $\cS$ be the set of weights which do not lie in $H$. Then $H=H_{\cS^c}$ and $\cS$ is minimal. Hence the assignment $\cS\mapsto H_{\cS^c}$ is a bijection between the minimal subsets of $\{1,..., n\}$ and the subspaces of $L^\vee_\R$ which are spanned by the weights they contain. 

It follows immediately from part (ii) of Definition/Lemma \ref{defn.relevant} that the subspace $H_{\cS^c}$ is relevant  iff there is a vector $k\in\Z^{\cS^c}$ with strictly positive entries which maps to zero under $Q$. Such a vector is exactly an element  $k\in N^\vee$ such that $k(A(e_i))=0$ if $i\in \cS$ and $k(A(e_i))>0$ if $i\in \cS^c$.  Since the polytope $\Pi$ lives in an affine hyperplane of height 1, the existence of such a $k$ is the statement that $\cS$ is all the rays on a face of $\Pi$. 
\end{proof}

The zero subspace $H={0}$ is always relevant, and since we assume there are no zero weights it corresponds to the whole polytope $\Pi$. The empty set is a face of $\Pi$ (in the sense of the above proof), it corresponds to the subspace $H=L^\vee_\R$ which is therefore relevant.

\subsection{The conjecture}

Let $W$ be a wall separating two chambers in the secondary fan. Recall that we have the following two objects associated to $W$:

\begin{enumerate}
\item A toric variety $Z_W$. The wall $W$ has an associated Higgs GIT problem, and $Z_W$ is the phase of this problem coming from a character on the relative interior of $W$. 

\item A toric rational curve $C_W$ in the secondary stack. $W$ is a codimension 1 cone in the secondary fan and $C_W$ is the associated curve.

\end{enumerate}

We can decompose $D^b(Z_W)$ using the algorithm of Section \ref{sec.SODs}, and the factors that appear are indexed by the relevant subspaces $H\subset L^\vee_\R$ contained in the hyperplane $\langle W \rangle$. Each such subspace defines a Higgs GIT problem with a non-empty minimal phase $Z_H$, and by Theorem \ref{thm.JH} the multiplicity of $D^b(Z_H)$ in $D^b(Z_W)$ is well-defined.
\pgap

Relevant subspaces correspond (by Proposition \ref{prop.minimaltorelevant}) to minimal faces $\Gamma$ of the polytope $\Pi$, and these in turn index the components of the discriminant locus.  As discussed in Section \ref{sec.multsintro} we are interested in the intersection of $\overline{\nabla}_\Gamma$ with the curve $C_W$.

\begin{lem}\label{lem.intersection} Let $\Gamma$ be a minimal face.  If $H_{\Gamma^c}$ is not contained in $W$ then $\overline{\nabla}_\Gamma$ doesn't meet the curve $C_W$.
\end{lem}
\begin{proof} Consider the projection map
$$\pi: L^\vee \to L_\Gamma^\vee$$
or its real version $L^\vee_\R \to (L_\Gamma^\vee)_\R$, whose kernel is $H_{\Gamma^c}$.  This map takes a stability condition for the original GIT problem and restricts it to give one for the Coloumb GIT problem associated to $\Gamma$. If we take a chamber of stability conditions and restrict them then they will all lie in a single chamber for the Coloumb GIT problem (if two stability conditions induce the same triangulation of $\Pi$ then they evidently induce the same triangulation of the face $\Gamma$). This says that $\pi$ is a map of fans, between the secondary fan for the original problem and the secondary fan for the Coloumb problem, hence it induces a toric morphism
$$\pi: \mfF \to \mfF_\Gamma$$
between the two secondary toric varieties.

Recall that $\nabla_\Gamma$ is defined as the preimage of the discriminant locus $\nabla'_\Gamma \subset (L_\Gamma)^\vee_{\C^*}$ under the projection $\pi: (L^\vee)_{\C^*} \to (L_\Gamma)^\vee_{\C^*}$.  Since $\pi$ extends to the toric boundary we can also say that $\overline{\nabla}_\Gamma \subset \mfF $ is contained in the preimage of $\overline{\nabla'}_{\Gamma} \subset \mfF_\Gamma $.

The wall $W$ is a codimension 1 cone in $L^\vee$. If it doesn't contain $H_{\Gamma^c}$ then $\pi(W)$ is a top-dimensional cone in $L_\Gamma^\vee$ so $\pi(C_W)$ is one of the toric fixed points in $\mfF_\Gamma$.  Corollary \ref{cor.missesfp} says that $\overline{\nabla'}_{\Gamma}$ avoids all the toric fixed points hence $\overline{\nabla}_\Gamma$ misses $C_W$. 
\end{proof}

\begin{conj}\label{conj.mainconj} Let $W\subset L^\vee_\R$ be a wall, let $\Gamma\subset \Pi$ be a minimal face and let $H=H_{\Gamma^c}$ the corresponding relevant subspace. Assume that $H\subseteq \langle W \rangle$. Write
\begin{enumerate}\setlength{\itemsep}{5pt}
\item[] $n_{\Gamma, W}$ for the multiplicity of $D^b(Z_H)$ in $D^b(Z_W)$, and
\item[] $m_{\Gamma, W}$ for the intersection multiplicity of $\overline{\nabla}_\Gamma$ with $C_W$.
\end{enumerate}
Then $n_{\Gamma, W} = m_{\Gamma, W}$.
\end{conj}

\begin{rem} We could allow the case when $H$ doesn't lie in $\langle W \rangle$: then $D^b(Z_H)$ is not a factor in $D^b(Z_W)$ so we should set $n_{\Gamma, W}=0$, and by Lemma \ref{lem.intersection} $m_{\Gamma, W}=0$ also. 
\end{rem}

We will now prove various special cases of this conjecture. The most straight-forward case is when $\rk L_\Gamma =1$ so $H$ is a hyperplane, hence $H=\langle W \rangle$. 

\begin{prop}\label{prop.H=W} If $\rk L_\Gamma =1$ then $n_{\Gamma, W} = m_{\Gamma, W} = 1$. 
\end{prop}
\begin{proof} In this case $Z_H$ is the minimal phase for the Higgs GIT problem that produces $Z_W$, so $n_{\Gamma, W}=1$. 

As in Lemma \ref{lem.intersection} we consider the map $\pi: \mfF \to \mfF_\Gamma$. This map induces an isomorphism from $C_W$ to $\mfF_\Gamma$. The discriminant locus $\nabla'_\Gamma\subset \mfF_{\Gamma}$ is a single non-fixed point (Example \ref{eg.rank1discriminant}) and $\overline{\nabla}_\Gamma$ is its pre-image, as a divisor on $\mfF$. So intersecting this divisor with $C_W$ gives $m_{\Gamma, W}=1$. 
\end{proof}

\begin{rem} This proposition includes the case when $L$ itself has rank 1 and hence $H=W$ is the origin. This is a  vacuous case of our conjecture:  $C_W$ is the whole of $\mfF$, there is only the principal component of $\overline{\nabla}$ which is a single point,  $Z_W$ is a point, and  the decomposition of $D^b(Z_W)$ is trivial.
\end{rem}

We can get a less trivial special case by increasing the rank by one.

\begin{prop} \label{prop.rank2} If $\rank L_\Gamma =2$ then  $n_{\Gamma, W} = m_{\Gamma, W}$. 
\end{prop}
\begin{proof}
In this case $H$ is a hyperplane in $\langle W \rangle$. Since the projection $\pi: L^\vee \to L^\vee_\Gamma$ is a map of fans $\pi(W)$ must be a ray, so $W$ must lie completely on one side of $H$. Pick a primitive one-parameter subgroup $\lambda$ normal to $W$, then the GIT problem producing $Z_W$ consists of the vector space $V^{\lambda}$ - these are weights that lie in $\langle W \rangle$ - acted on by the torus $T/\lambda$. Then $H$ is normal to some primitive one-parameter subgroup $\mu \in L/\langle\lambda\rangle$ and we orient $\mu$ so that it pairs positively with $W$. Recall that $Z_W$ is defined to be the phase associated to a generic character in $W$, so such a character pairs positively with $\mu$. For the decomposition of $D^b(Z_W)$ the important quantity is
$$\kappa \;=\; (\det V^\lambda)(\mu) \;=\;  \sum_{\mathrm{weights}\, q_i \in \langle W \rangle } q_i(\tilde{\mu}) $$
for any $\tilde{\mu}\in L$ lifting $\mu$, and 
$$n_{\Gamma, W} = \max\{ \kappa, 0\}$$
(see Remark \ref{rem.mult}). 

Now we compute the intersection multiplicity $m_{\Gamma, W}$. To start with, let's assume that $\Gamma=\Pi$ so $L$ itself has rank 2 and $H$ is the origin. Then we wish to compute the intersection multiplicity of the principal component $\overline{\nabla}_A$ with the boundary curve $C_W$. To do this we use the Horn uniformization map
$$\P(L_\Gamma)_\C \dashrightarrow \overline{\nabla}_A $$
from Section \ref{sec.Horn}, which in this case is actually a morphism since $\P(L_\C) \cong \P^1$. In explicit co-ordinates, as a rational map to $(\C^*)^2$, this is given by:
$$\lambda_1: \lambda_2  \;\mapsto \; \left( \prod_{i=1}^n \Big(  Q_{i1}\lambda_1 + Q_{i2}\lambda_2\Big)^{Q_{i1}}\!\!, \; \;  \prod_{i=1}^n \Big(  Q_{i1}\lambda_1 + Q_{i2}\lambda_2\Big)^{Q_{i2}}      \right) $$
Without loss of generality we may assume that $W$ is the ray through $(1,0)$. This ray in the secondary stack corresponds to a partial compactification of the torus $(\C^*)^2$, the subset:
$$\C \times \C^* \subset \mfF$$
The subset where the first co-ordinate is zero is $C_W$ with its fixed points deleted.  Since  $\overline{\nabla}_A$ avoids the fixed points the only way that $\lambda_1\!:\!\lambda_2$ can map to $C_W$ is if there exists an $i$ such that $Q_{i1}\lambda_1 + Q_{i1} \lambda_2 = 0$  and $Q_{i2}=0$, hence $\lambda_1\!:\!\lambda_2=0\!:\!1$.   Then the intersection multiplicity is given by
$$\sum_{i|\, Q_{i2}=0} Q_{i1}$$
if this sum is strictly positive, and zero otherwise. But these rows of $Q$ are precisely the weights $q_i$ that lie on $\langle W \rangle$, and we may set $\tilde{\mu}=(1,0)^\top$, so this sum is $\kappa$ and hence $m_{\Gamma, W}=n_{\Gamma, W}$ in this case.

To finish we must compute $m_{\Gamma, W}$ for $\rank L_{\Gamma}=2$ but $\Gamma \varsubsetneq \Pi$.  Once again we use the projection $\pi: \mfF \to \mfF_{\Gamma}$. It maps $W$ to a wall $W'$ for the rank 2 GIT problem, and hence it maps $C_W$ isomorphically (at least away from its fixed points) onto the boundary curve $C_{W'}  \subset \mfF_{\Gamma}$. Since $\overline{\nabla}_\Gamma$ is the preimage of $\overline{\nabla'}_{\Gamma}\subset \mfF_{\Gamma}$  it is enough to compute the intersection multiplicity of $\overline{\nabla}_{\Gamma}$ with $C_{W'}$ inside the two-dimensional space $\mfF_{\Gamma}$. But this was the calculation we just performed, and the result is still $\max\{\kappa, 0\}$ since $q_i\in \langle W\rangle$ iff $\pi(q_i)\in \langle W' \rangle$.
\end{proof}

\begin{rem}The above result and its proof are quite close to \cite[Prop. 4.4.]{HLSh}.\end{rem}

\begin{thm}\label{thm.rank2} If $\rank L =2$ then Conjecture \ref{conj.mainconj} holds. 
\end{thm}

\begin{proof} The wall $W$ is a ray and there are only two possibilities for $H$: either $H=\langle W \rangle$ (if this is a relevant subspace) or $H={0}$. The first case is covered by Proposition \ref{prop.H=W} and the second by Proposition \ref{prop.rank2}. 
\end{proof}

The main obstacle to extending our proofs to higher rank is the fact that Horn uniformization may no longer be a morphism so it becomes harder to compute the intersection multiplicity $m_{\Gamma, W}$. However in special cases it is still possible to verify the conjecture - see \cite[Sect.~10.2]{Kite} for some more examples.

% ----------------------------------------------------------------
\bibliographystyle{halphanum}

\end{document}